\documentclass{amsart}

\usepackage{amssymb}
\usepackage{enumerate,amsfonts,calc,amsmath,verbatim}

\usepackage[all]{xy}
\SelectTips{eu}{}

\newcommand{\Ext}{\operatorname{Ext}}

\newcommand{\Tor}{\operatorname{Tor}}
\newcommand{\tor}[4][{}]{\operatorname{Tor}_{#1}^{#2}(#3,#4)}
\newcommand{\gtor}[5]{\operatorname{Tor}_{#1}^{#2}(#3,#4)_{#5}}
\newcommand{\Hom}{\operatorname{Hom}}

\newcommand{\Coker}{\operatorname{Coker}}

\newcommand{\rank}{\operatorname{rank}}

\newcommand{\agr}[2][{}]{{{#2}^{\mathsf g}_{#1}}}

\newcommand{\lin}{\operatorname{lin}}

\newcommand{\shift}{\mathsf\Sigma}

\newcommand{\ov}{\overline}

\newcommand{\col}{\colon}

\newcommand{\ges}{\geqslant}
\newcommand{\les}{\leqslant}
\newcommand{\Ker}{\operatorname{Ker}}

\newcommand{\BZ}{{\mathbb Z}}

\newcommand{\fm}{{\mathfrak m}}

\newcommand{\fp}{{\mathfrak p}}
\newcommand{\fq}{{\mathfrak q}}
\newcommand{\fs}{{\mathfrak s}}
\newcommand{\ft}{{\mathfrak t}}

\newcommand{\xra}{\xrightarrow}
\newcommand{\lra}{\longrightarrow}

\newcommand{\bm}{\begin{matrix}}
\newcommand{\dm}{\end{matrix}}

\theoremstyle{remark}

\theoremstyle{plain}
\swapnumbers \theoremstyle{plain}

\newtheorem{theorem}{Theorem}[section]

\newtheorem{proposition}[theorem]{Proposition}
\newtheorem{lemma}[theorem]{Lemma}

\newtheorem{corollary}[theorem]{Corollary}

\theoremstyle{definition}

\newtheorem{chunk}[theorem]{}

\newtheorem{example}[theorem]{Example}

\theoremstyle{remark}

\newtheorem{remark}[theorem]{Remark}

\newtheorem{construction}[theorem]{Construction}

\numberwithin{equation}{theorem}
\numberwithin{subchunk}{theorem}

\newcommand{\dd}{\partial}

\newcommand{\hh}[2][{}]{{\operatorname{H}_{#1}(#2)}}

\newcommand{\hilb}[2]{H_{#1}(#2)}
\newcommand{\po}[3][R]{P^{#1}_{#2}(#3)}

\newcommand{\zz}{{\operatorname{Z}}}

\newcommand{\bsh}{{\boldsymbol h}}

\begin{document}

\title[Free resolutions over short local rings]{Free resolutions over  
short local
rings}

   \begin{abstract}
     The structure of minimal free resolutions of finite modules $M$ over commutative
     local rings $(R,{\mathfrak m},k)$ with ${\mathfrak m}^3=0$ and
     $\operatorname{rank}_k({\mathfrak m}^2)<\operatorname{rank}_k ({\mathfrak
       m}/{\mathfrak m}^2)$ is studied.  It is proved that over generic $R$ every $M$ has
     a Koszul syzygy module.  Explicit families of Koszul modules are identified.  When
     $R$ is Gorenstein the non-Koszul modules are classified.  Structure theorems are
     established for the graded $k$-algebra $\operatorname{Ext}_R(k,k)$ and its graded
     module $\operatorname{Ext}_R(M,k)$.
   \end{abstract}

\thanks{Research partly supported by NSF grants DMS 0201904 (LLA) and  
DMS 0602498
(SI)}

\author[L.~L.~Avramov]{Luchezar L.~Avramov}
\address{Luchezar L.~Avramov\\ Department of Mathematics\\
   University of Nebraska\\ Lincoln\\ NE 68588\\ U.S.A.}
     \email{avramov@math.unl.edu}

\author[S.~B.~Iyengar]{Srikanth B.~Iyengar}
\address{Srikanth B.~Iyengar\\ Department of Mathematics\\
   University of Nebraska\\ Lincoln\\ NE 68588\\ U.S.A.}
     \email{iyengar@math.unl.edu}

\author[L.~M.~\c{S}ega]{Liana M.~\c{S}ega}
\address{Liana M.~\c{S}ega\\ Department of Mathematics and Statistics\\
   University of Missouri\\ \linebreak Kansas City\\ MO 64110\\ U.S.A.}
     \email{segal@umkc.edu}


\dedicatory{To the memory of our friend and colleague Anders Frankild.}

\maketitle

\section*{Introduction}

This paper is concerned with the structure of minimal free resolutions
of finite (that is, finitely generated) modules over a commutative 
noetherian local ring $R$ whose maximal ideal $\fm$ satisfies $\fm^3=0$.  
Over the last 30 years this special class has emerged as a testing ground for properties
of infinite free resolutions.  For a finite module $M$ such properties are 
often stated in terms of its \emph{Betti numbers} $\beta_n^R(M)=
\rank_k\Ext^n_R(M,k)$, where $k=R/\fm$, or in terms of its
\emph{Poincar\'e series}
  \[
\po Mt=\sum_{i=0}^\infty\beta_n^R(M)\,t^n\in\BZ[\![t]\!]\,.
  \]

  Patterns that had been conjectured not to exist at one time have
  subsequently been discovered over rings with $\fm^3=0$: Not finitely generated algebras
  $\Ext_R(k,k)$ (Roos, 1979); transcendental Poincar\'e series $\po kt$ (Anick, 1980);
  modules with constant Betti numbers and aperiodic minimal free resolutions (Gasharov and
  Peeva, 1990); families of modules with rational Poincar\'e series that admit no common
  denominator (Roos, 2005); reflexive modules $M$ with
  $\Ext_R^{n}(M,R)=0\ne\Ext^{n}_R(\Hom_R(M,R),R)$ for all $n\ge1$ (Jorgensen and \c Sega,
  2006).

On the other hand, important conjectures on infinite free resolutions that are still open
in general have been verified over rings with $\fm^3=0$: Each sequence
$(\beta^R_n(M))_{n\ges0}$ is eventually non-decreasing, and grows either polynomially or
exponentially (Lescot, 1985). When $M$ has infinite projective dimension one has
$\Ext_R^{n}(M,M\oplus R)\ne0$ for infinitely many $n$ (Huneke, \c Sega and Vraciu, 2004).

The work presented below is motivated by an `unusually high' incidence in the appearance of modules $M$ with `Koszul-like' behavior, exemplified by an equality
   \[
\po Mt=\frac{p_M(t)}{1-et+rt^2} \tag{$\ast$}
   \]
with $e=\rank_k(\fm/\fm^2)$, $r=\rank_k(\fm^2)$, and $p_M(t)\in\BZ[t]$.

We provide structural explanations for the phenomenon and natural 
conditions for its occurrence.  To describe our results, note that the 
property $\dd(F)\subseteq \fm F$ of a minimal free resolution $F$ of $M$ 
allows one to form for each $j\ge0$ a complex
   \[
\lin_j(F)=\quad0\to
\frac{F_{j}}{\fm F_{j}}\to
\cdots\to
\frac{\fm^{j-n}F_{n}}{\fm^{j+1-n}F_{n}}\to
\cdots\to
\frac{\fm^{j}F_{0}}{\fm^{j+1}F_{0}}\to
\frac{\fm^{j}M}{\fm^{j+1}M}\to
0
\]
of $k$-vector spaces.   Following \cite{HI}, we say that $M$ is
\emph{Koszul} if every complex $\lin_j(F)$ is acyclic.  When $R$ is a  
graded $k$-algebra generated in degree $1$ and $M$ is a graded 
$R$-module generated in a single degree, say $d$, this means that 
$M$ has a $d$-linear free resolution.  The results in this paper are 
new also in this more restrictive setup.

Our main theorem states that if there exists an $x\in\fm$ satisfying $x^2=0$ and
$\fm^2=x\fm$, then every finite $R$-module has a syzygy module that is
Koszul; see Section \ref{sec:Minimal resolutions}.  The crucial step
is to  find a quadratic hypersurface ring mapping onto $R$ by a Golod
homomorphism.  Results of Herzog and Iyengar \cite{HI} then apply and
yield, in particular,  formula ($\ast$).

We call an element $x$ as above a \emph{Conca generator} of $\fm$  
because Conca \cite{Co1} shows that such an $x$ exists in generic 
standard  graded $k$-algebras with $r\le e-1$.

Assuming that $\fm$ has a Conca generator, in Section \ref{sec:Structure of cohomology} we provide a finite presentation of the algebra $\mathcal E=\Ext_R(k,k)$ and prove that for every finite $R$-module $M$ the $\mathcal E$-module $\Ext_R(M,k)$ has a resolution of length at most $2$
by finite free graded $\mathcal E$-modules.  This yields information
on the degree of the polynomial $p_M(t)$ in ($\ast$).

The last two sections are taken up by searches for Koszul modules.

In Section \ref{sec:Koszul modules} we prove that modules
annihilated by a Conca generator are Koszul.

In Section \ref{sec:Gorenstein rings and related rings}, we identify 
the rings with $\fm^3=0$ and $\rank_k\fm^2\le1$ whose maximal 
ideal $\fm$ has a Conca generator.  The Gorenstein rings $R$ with $e\ge2$
are among them; in particular, they satisfy formula $(*)$, which is
known from Sj\"odin \cite{Sj2}.  We prove that over such rings
the indecomposable non-Koszul modules are precisely 
the negative syzygies of $k$, and that these modules are characterized by 
their Hilbert function.  For $e=2$ this follows from Kronecker's
classification of pairs of commuting matrices.  It is unexpected that 
the classification extends \emph{verbatim} to Gorenstein rings with 
$\fm^3=0$ and $e\ge3$, which have wild representation type.

\section{Minimal resolutions}
\label{sec:Minimal resolutions}

In this paper the expression \emph{local ring} $(R,\fm,k)$ refers to
a commutative Noetherian ring $R$ with unique maximal ideal $\fm$ and
residue field $k=R/\fm$.  Such a ring $R$ is called \emph{Koszul} if its
residue field is a Koszul $R$-module, as defined in the introduction;
when $R$ is standard graded this coincides with the classical notion.
An \emph{inflation} $R\to R'$ is a homomorphism to a local ring
$(R',\fm',k')$ that makes $R'$ into a flat $R$-module and satisfies
$\fm'=R'\fm$.  

Here is the main result of  this  section:
\begin{theorem}
  \label{thm:koszuls}
Let $(R,\fm,k)$ be a local ring.

If for some inflation $R\to(R',\fm',k')$ the ideal $\fm'$ has a
Conca generator, then $R$ is Koszul, each finite $R$-module $M$
has a Koszul syzygy module, and one has
   \begin{align}
  \label{eq:koszul-field}
\po kt=&\frac{1}{1-et+rt^2}\,;
  \\ \label{eq:koszul-module}
\po Mt=&\frac{p_M(t)}{1-et+rt^2}
  \quad\text{with}\quad
p_M(t)\in\BZ[t] \,,
   \end{align}
where we have set $e=\rank_k(\fm/\fm^2)$ and $r=\rank_k(\fm^2/\fm^3)$.
      \end{theorem}
A notion of Conca generator was introduced above;
we slightly expand it:

\begin{chunk}
   \label{conca}
We say that $x$ is a \emph{Conca generator} of an ideal $J$ if $x$
is in $J$ and satisfies $x\ne0=x^2$ and $xJ=J^2$.   One then has
$J^3=(xJ)J=x^2J=0$, and also $x\notin J^2$: The contrary would imply
$x\in J^2=xJ\subseteq J^3=0$, a contradiction.
  \end{chunk}

When the ring $R$ in the theorem is standard graded a result of Conca,
Rossi, and Valla, see \cite[(2.7)]{CRV}, implies that it is Koszul;
the proof relies on the theory of Koszul filtrations developed by
these authors.  A direct proof is given by Conca \cite[Lem.~2]{Co1}.
Neither argument covers the local situation, nor gives information on
homological properties of $R$-modules other than $k$.

We prove Theorem~\ref{thm:koszuls} via a result on the structure of $R$.
It is stated in terms of Golod homomorphisms, for which we recall one
of several possible definitions.

\begin{chunk}
   \label{massey}
Let $(Q,\fq,k)$ be a local ring and $D$ a minimal free resolution of $k$
over $Q$ that has a structure of graded-commutative DG algebra; one 
always exists, see \cite{Gu1}.

Let $\varkappa\col Q\to R$ be a surjective homomorphism of rings, and
set $A=R\otimes_QD$.  For $a\in A$ let $|a|=n$ indicate $a\in A_n$,
and set $\ov a=(-1)^{|a|+1}a$.  Let $\bsh$ denote a homogeneous
basis of the graded $k$-vector space $\hh[\ges1]A$.

The homomorphism $\varkappa$ is \emph{Golod} if there is a function
$\mu\col\bigsqcup_{m=1}^{\infty}\bsh^m\to A$ satisfying
\begin{align}
\label{eq:tmo1}
\mu&(h)
\text{ is a cycle in the homology class of $h$
for each } h\in\bsh\,;
   \\
\label{eq:tmo2}
\dd\mu&(h_{1},\dots,h_{m})
=\sum_{i=1}^{m-1}\ov{\mu(h_{1},\dots,h_{i})}
\mu(h_{i+1},\dots,h_{{m}})
\quad\text{for each }m\geq 2\,;
    \\
   \label{eq:tmo3}
\mu&(\bsh^m)\subseteq \fm A\quad\text{for each }m\geq 1\,.
\end{align}
   \end{chunk}

The following structure theorem is used multiple times in the paper.

   \begin{theorem}
     \label{thm:factorization}
Let $(R,\fm,k)$ be a local ring and $x$  a Conca generator of $\fm$.

There  exist a regular local ring $(P,\fp,k)$, a minimal set of
generators $u_1,\dots,u_e$ of $\fp$, and a surjective homomorphism
$\pi\col P\to R$, with $\Ker(\pi)\subseteq\fp^2$ and  $\pi(u_e)=x$.

Any such $\pi$ induces a Golod homomorphism $Q\to R$, where $Q=P/ 
(u_e^2)$.
   \end{theorem}

Indeed, $R$ is artinian by \ref{conca}, so Cohen's Structure Theorem
yields a surjection $\pi\col P\to R$  with  $(P,\fp,k)$ regular local  
and $\Ker(\pi)\subseteq\fp^2$; now \ref{conca} gives $x\notin\fm^2$, 
hence any $u_e\in\fp$ with $\pi(u_e)=x$ extends to a minimal  
generating set of $\fp$.  It remains to prove that $Q\to R$ is Golod.  For
this we use a construction going back to H.~Cartan.

    \begin{construction}
    \label{con:resolution}
Form a complex of free $Q$-modules of rank one
   \[
D=\quad \cdots\lra Qy_{n}\xra{\ u_e\ }Qy_{n-1}\lra\cdots \lra
Qy_1\xra{\ u_e\ }Qy_0\lra 0\lra\cdots
   \]
where $y_n$ is a basis element in degree $n$. The multiplication table
\begin{gather*}
y_{2i}y_{2j}=\binom{i+j}{i}y_{2(i+j)}=y_{2j}y_{2i}\\
y_{2i}y_{2j+1}=\binom{i+j}{i}y_{2(i+j)+1}=y_{2j+1}y_{2i}\\
y_{2i+1}y_{2j+1}=0=y_{2j+1}y_{2i+1}
\end{gather*}
for all $i,j\ge0$ turns $D$ into a graded-commutative
DG $Q$-algebra with unit $y_0$.
      \end{construction}

\begin{lemma}
\label{lem:resolution}
Let $C$ be the Koszul complex on the image in $Q$ of
$\{u_1,\dots,u_{e-1}\}$.

The DG $Q$-algebra $C\otimes_QD$ is a minimal free resolution of $k$.

The DG $R$-algebra $A=B\otimes_QD$, where $B=R\otimes_QC$,
satisfies the conditions:
    \[
\dd(A_1)=\fm\subseteq R=A_0\,,\quad \fm^2A_{\ges1}\subseteq
\dd(\fm A_{\ges 2})\,,\quad\text{and}\quad \zz_{\ges
1}(A)\subseteq\fm A_{\ges 1}\,.
    \]
    \end{lemma}

    \begin{proof}
The canonical augmentation $D\to Q/u_eQ$ is a quasi-isomorphism.
It induces the first quasi-isomorphism below, because $C$ is a finite
complex of free $Q$-modules:
    \[
C\otimes_QD\simeq C\otimes_Q(Q/u_eQ)\simeq k\,.
    \]
The second one holds because the sequence $u_1,\dots,u_{e-1}$ is
$(Q/u_eQ)$-regular.

The relations $\dd(A_1)=\fm\subseteq R=A_0$ hold by construction.

Every element $a\in A_n$ can be written uniquely in the form $a=\sum_{i=0}^n b_i\otimes y_i$ with
$b_i\in B_{n-i}$.  Supposing that it is a cycle, one gets
   \[
\dd(b_i)=(-1)^{n-i}x b_{i+1}
  \quad\text{for}\quad
i=0,\dots,n-1
    \quad\text{and}\quad
  \dd(b_n)=0\,.
   \]
If $b\in B$ satisfies $\dd(b)\in\fm^2B$, then $b$ is in $\fm B$; see \cite[Cor~2.1]{Lo}.  Thus, $\dd(b_n)=0$ implies $b_n\in\fm B$. Assuming by descending induction that $b_j\in \fm B$ for some $j\le n$, the equalities above imply $\dd(b_{j-1})\in\fm^2B$, whence $b_{j-1}\in\fm B$,
and finally $a\in\fm A$.

It remains to verify the condition $\fm^2A\subseteq \dd(\fm A)$. One has
$B_{\les {e-1}}=B$, so to obtain the desired inclusion we induce on $n$
to show that the following holds:
    \[
\fm^2(B_n\otimes_QD)\subseteq \dd(\fm(B_{\les n}\otimes_QD))
\quad\text{for each } n\,.
    \]
For $n\le-1$ the formula above holds because $B_{<0}=0$, so assume it
holds for some integer $n\ge-1$.  As $B$ and $D$ are complexes of free
modules, one gets  relations:
    \begin{align*}
      x(B_{n+1}\otimes_QD) &= B_{n+1}\otimes_Q u_eD
\\
    &\subseteq B_{n+1}\otimes_Q\dd(D)
\\
    &\subseteq \dd(B_{n+1})\otimes_QD + \dd(B_{n+1}\otimes_QD)
\\
    &\subseteq \fm (B_{n}\otimes_QD) + \dd(B_{n+1}\otimes_QD)\,.
    \end{align*}
The resulting inclusion provides the second link in the following chain:
    \begin{align*}
      \fm^2(B_{n+1}\otimes_QD) &= \fm x(B_{n+1}\otimes_QD)
\\
     &\subseteq \fm^2(B_{n}\otimes_QD) + \fm \dd(B_{n+1}\otimes_QD)
\\
     &\subseteq \dd(\fm (B_{\les n}\otimes_QD))+
     \dd(\fm(B_{n+1}\otimes_QD))
\\
      &=\dd(\fm(B_{\les n+1}\otimes_QD))\,.
    \end{align*}
The induction hypothesis gives the third link, so the argument is
complete.
    \end{proof}

\begin{proof}[Proof of Theorem \emph{\ref{thm:factorization}}]
We construct a function $\mu$ as in \ref{massey} by induction on $m$.

To start, choose for each $h\in\bsh$ a cycle $z_h$ in the homology
class of $h$; the last inclusion in Lemma \ref{lem:resolution} yields
$z_h\in\fm A$.  Suppose that $\mu$ has been defined on $\bsh^{m-1}$ for
some integer $m\geq 2$.  The induction hypothesis then provides an  
element
   \[
\sum_{i=1}^{m-1}\ov{\mu(h_{1},\dots,h_{i})}
    \mu(h_{i+1},\dots,h_{m})\in\fm^2A\,.
   \]
As $\fm^2A\subseteq \dd(\fm A)$ holds, by Lemma~\ref{lem:resolution},
the element above is the boundary of some element in $\fm A$, which we
name $\mu(h_{1},\dots,h_{m})$.  The induction step is complete.
    \end{proof}

We will use an adaptation of a construction of Golod proposed by
Gulliksen:

\begin{construction}
  \label{golod-construction}
For each integer $n\ge1$ set $\bsh_n=\{h\in\bsh : |h|=n\}$, and choose a
free $R$-module $V_{n}$ with basis
$\{v_h\}_{h\in\bsh_n}$.  The free $R$-modules
    \begin{equation}
   \label{eq:golod-basis}
G_n=\bigoplus_{l+m+n_1+\dots+n_m=n} A_l\otimes_R
V_{n_1}\otimes_R\cdots\otimes_R V_{n_m}
    \end{equation}
and the $R$-linear maps $\dd\col G_n\to G_{n-1}$ defined by the formula
      \begin{equation}
    \label{eq:golod-differential}
      \begin{split}
\dd(a\otimes v_{h_1}\otimes\cdots\otimes v_{h_m}) = &\
\dd(a)\otimes v_{h_1}\otimes\cdots\otimes v_{h_m} \\
    & +(-1)^{|a|}\sum_{i=1}^m
      a\mu(h_{1},\dots,h_{i})\otimes v_{h_{i+1}}\otimes
      \cdots\otimes v_{h_m}
    \end{split}
      \end{equation}
then form a minimal free resolution of $k$ over $R$; see
\cite[Prop.~1]{Gu2}.
   \end{construction}

For a finite $R$-module $M$ and $n\ge0$ set $\Omega^R_n(M)=
\Coker(\dd^F_{n+1})$, where $F$ is a minimal free resolution of $M$;
up to isomorphism, $\Omega^R_n(M)$ does not depend on $F$.

  \begin{remark}
   \label{rem:inflation}
Let $R\to(R',\fm',k')$ be an inflation.

Let $M$ be a finite $R$-module and set $M'=M\otimes_RR'$.

For each $n\ge0$ one has $\fm'^nM'/\fm'^{n+1}M'\cong
k'\otimes_k(\fm^nM/\fm^{n+1}M)$.

If $F$ is a minimal free resolution of $M$ over $R$, then $R'\otimes_RF$
is one of $M'$ over  $R'$, so $\po Mt=\po[R']{M'}t$ and
$\Omega^{R'}_n(M')\cong R'\otimes_R\Omega^{R}_n(M)$.  The complexes
$\lin_j(F)$ from the introduction satisfy $\lin^{R'}_j(R'\otimes_{R}F)\cong
k'\otimes_k \lin^{R}_j(F)$ for each $j$.  Thus, $\lin_j(R'\otimes_{R}F)$ and
$\lin_j(F)$ are exact simultaneously, so $M$ is Koszul if and only if
$M'$ is.
  \end{remark}

   \begin{remark}
   \label{rem:syzygy}
Set $\hilb Mt=\sum_{n=0}^\infty\rank_k(\fm^nM/\fm^{n+1}M)\cdot t^n$.

When $R$ is a Koszul ring with $\fm^{l}=0$ for some integer $l$, and $N=\Omega^R_i(M)$ is a Koszul module
the proof of \cite[(1.8)]{HI} shows that $\po Mt\cdot\hilb R{-t}-t^{i}\hilb N{-t}$ is a 
polynomial of degree at most $i+l-2$.
   \end{remark}

   \begin{proof}[Proof of Theorem \emph{\ref{thm:koszuls}}]
In view of Remark \ref{rem:inflation}, we may replace $R$ with $R'$
and $M$ with $M'$, and so assume that the ideal $\fm$ itself has a
Conca generator.

We start by proving that the ring $R$ is Koszul.  By \cite[(7.5)]{Se},
it suffices to show that for all integers $n\ge1$ and $j\ge1$ the map
    \[
\Tor_n^R(\iota_j,k)\col \Tor_n^R(\fm^{j+1},k)\to\Tor_n^R(\fm^{j},k)\,,
   \]
induced by the inclusion $\iota_j\col\fm^{j+1}\to\fm^j$ is equal to $0$.
Because we have $\fm^3=0$, see \ref{conca}, this boils down to proving
$\Tor_n^R(\iota_1,k)=0$ for all $n\ge1$.  These maps are the  
homomorphisms
induced in homology by the inclusion $\fm^2 G\subseteq\fm G$, where $G$
is a free resolution of $k$ over $R$.  We choose the one from  
Construction
\ref{golod-construction}.

As seen from \eqref{eq:golod-basis}, the graded $R$-module $\fm^2G$
is spanned by elements
    \[
g=r a\otimes v_{h_1}\otimes \cdots \otimes v_{h_m}
    \]
with $r\in \fm^2$, $a\in A$, and $h_i\in \bsh$. By induction on $m$ we
prove $g\in\dd(\fm G)$.  Lemma \ref{lem:resolution} settles the case
$m=0$.  Suppose $m\ge1$. The base case yields $ra=\dd(b)$ for some $b\in\fm
A$, so from \eqref{eq:golod-differential} one obtains:
   \begin{gather*}
g= \dd(b) \otimes v_{h_1}\otimes\cdots \otimes v_{h_m}
       = \dd(b\otimes v_{h_1}\otimes\cdots \otimes v_{h_m})-(-1)^{|b|} 
c\,,
   \\
\text{where}\quad c=\sum_{i=1}^mb\mu(h_{1},\dots,h_{i})\otimes
v_{h_{i+1}}\otimes \cdots\otimes v_{h_m}\,.
     \end{gather*}
One has $\mu(h_{1},\dots,h_{i})\in\fm A$, see \eqref{eq:tmo3}, so $c$
is in $\fm^2G$. The induction hypothesis now yields $c\in\dd(\fm G)$,
so one gets $g\in\dd(\fm G)$, and the inductive proof is complete.

Theorem \ref{thm:factorization} yields a hypersurface ring $Q$ and a  
Golod
homomorphism $Q\to R$, so $M$ has a Koszul syzygy by \cite[(5.9)]{HI}.
Formulas \eqref{eq:koszul-module} and \eqref{eq:koszul-field} now follow
from Remark \ref{rem:syzygy} because one has $\fm^n=0$ for $n\ge3$;
see \ref{conca}.
    \end{proof}

   \begin{remark}
   \label{name}
Let $k$ be an algebraically closed field.  The quotient algebras of
the polynomial ring $k[X_1,\dots,X_e]$ by $\binom{e+1}2-r$ linearly
independent quadratic forms are parametrized by the points of an
appropriate Grassmannian. Conca proves that it contains a non-empty
Zariski open subset, whose points correspond to algebras $R$ for  
which the ideal $(X_1,\dots,X_e)R$ has a Conca generator; see 
\cite[Thm.~10]{Co1}.
   \end{remark}

A standard graded ring $R=k[X_1,\dots,X_e]/I$ is said to be
\emph{G-quadratic} if the ideal $I$ is generated by a Gr\"obner basis of
quadrics with respect to some system of coordinates and some term order.
The following implications hold, see \cite[Lem.~2]{Co1} for a proof of
the first one and \cite[2.2]{BHV} for one of the second:
   \[
(X_1,\dots,X_e)R\text{ has a Conca generator } \implies R\text{ is
G-quadratic } \implies R\text{ is Koszul}\,.
   \]

To show that the conclusions of Theorem \ref{thm:koszuls} may fail
over G-quadratic rings, even when $r\le e-1$ holds, we use a general
change-of-rings formula:

    \begin{remark}
   \label{rem:dress}
When $R$ is the fiber product  $S\times_kT$ of local rings $(S,\fs,k)$
and $(T,\ft,k)$ Dress and Kr\"amer \cite[Thm.~1]{DK} show that each
finite $S$-module $M$ satisfies
   \[
\po Mt=\frac{\po[S]Mt\cdot\po[T]kt}{\po[S]kt+\po[T]kt-\po[S]kt\po[T] 
kt}\,.
   \] \end{remark}

  \begin{example}
   \label{ex:non-koszul}
For a field $k$, an integer $e\ge5$, and indeterminates $X_1,\dots,X_e$ set
   \[
R=\frac{k[X_1,\dots,X_e]}{(X_1,X_2)^2+(X_3,X_4)^2
+(X_5,\dots,X_e)(X_1,\dots,X_e)}\,.
   \]

This $k$-algebra is evidently G-quadratic, local and has $\hilb Rt=1+et+4t^2$.
It is isomorphic to the fiber product  $S\times_kT$ of the
Koszul algebras
   \[
S=\frac{k[X_1,\dots,X_4]}{(X_1,X_2)^2+(X_3,X_4)^2}
   \quad\text{and}\quad
T=\frac{k[X_5,\dots,X_e]}{(X_5,\dots,X_e)^2}\,.
   \]
Thus, $\po[S]kt=(1-4t+4t^2)^{-1}$ and $\po[T]kt=(1-(e-4)t)^{-1}$, so
Remark \ref{rem:dress} yields
   \[
\po Mt=\po[S]Mt\cdot\frac{1-4t+4t^2}{1-et+4t^2}\,.
   \]
Roos \cite[\S5]{Ro} produces an infinite family of $S$-modules whose
Poincar\'e series admit no common denominator, so the conclusions of
Theorem \ref{thm:koszuls} fail for $R$.
   \end{example}

\section{Structure of cohomology}
\label{sec:Structure of cohomology}

In this section we obtain precise and explicit information on the graded
$k$-algebra $\Ext_R(k,k)$ and its graded left module $\Ext_R(M,k)$,
assuming that $\fm$ has a Conca generator.  We start by recalling
terminology and introducing notation.

Let $G$ be a free resolution of $k$ over $R$.  We identify $\Ext_R(k,k)$
with the cohomology, $\mathcal E$, of the complex $\Hom_R(G,G)$.
Composition gives $\Hom_R(G,G)$ a structure of a DG algebra, so $ 
\mathcal
E$ inherits a structure of graded algebra over $\mathcal E^0=k$.

  \begin{chunk}
   \label{rem:quadrics}
Let $(R,\fm,k)$ be a local ring and $x$ a Conca generator of $\fm$.

Multiplication with $x$ induces a surjection
$\upsilon\col\fm/\fm^2\to\fm^2$ with $\upsilon(x)=0$.  Thus, one
can choose $x_{1},\dots,x_{e}\in\fm\smallsetminus\fm^2$ so that
$x_{1}x,\dots,x_{r}x$ form a $k$-basis of $\fm^2$, the images of
$x_{r+1},\dots,x_{e}$ in $\fm/\fm^2$ form one of $\Ker(\upsilon)$, and
$x_e=x$.  It follows that $x_{1},\dots,x_{e}$ minimally generate $\fm$
and uniquely define elements $a^{ij}_h\in k$ satisfying
   \begin{equation}
  \label{eq:quadrics-local}
x_ix_j=\sum _{h=1}^r a^{ij}_hx_hx_e
  \quad\text{for}\quad 1\le i\le j\le e-1\,.
   \end{equation}

Let $\xi_1,\dots, \xi_e$ denote the basis of $\mathcal E^1$ dual to the
basis $x_{1},\dots,x_{e}$ of $\fm/\fm^2$, under the canonical  
isomorphism
$\mathcal E^1 \cong \Hom_k(\fm/\fm^2,k)$.
  \end{chunk}

We let $k\{\xi_1,\dots, \xi_e\}$ denote the non-commutative polynomial
ring on indeterminates $\xi_1,\dots,\xi_e$, set $[\xi_i,\xi_j]=
\xi_i\xi_j+\xi_j\xi_i$, and write $(\phi_1,\dots,\phi_r)$ for
the two-sided ideal generated by non-commutative polynomials
$\phi_1,\dots,\phi_r$.   A graded module $\mathcal M$ is said to be
\emph{bounded below} if $\mathcal M^n=0$ holds for all $n\ll0$.

In the proof of the next theorem, and later on, we use
non-traditional notation for associated graded objects:   For a local 
ring $(R,\fm,k)$ and an $R$-module $M$, set
\[
\agr R =
\bigoplus_{j\ges 0}\fm^{j}/\fm^{j+1}
   \quad\text{and}\quad
\agr M=
   \bigoplus_{j\ges 0}\fm^{j}M/\fm^{j+1}M\,.
\]
Of course, $\agr R $ is a graded ring and $\agr M$ is a graded module over it.

\begin{theorem}
  \label{thm:distinguished}
With the hypotheses and notation above the following hold.
  \begin{enumerate}[\rm(1)]
   \item
There is an isomorphism of graded $k$-algebras
   \begin{gather*}
\mathcal E\cong\frac{k\{\xi_1,\dots, \xi_e\}}{(\phi_1,\dots,\phi_r)}\,,
   \quad\text{where}\quad
  \\
\phi_h=[\xi_h,\xi_e]+\sum_{1\le i< j\le e-1}a^{ij}_h[\xi_i,\xi_j]
+\sum_{i=1}^{e-1}a^{ii}_h\xi_i^2
   \quad\text{for}\quad  h=1,\dots,r\,,
   \end{gather*}
and $a^{ij}_h$ are the elements in $k$ defined by formulas
\eqref{eq:quadrics-local}.
   \item
The images of the words in $\{\xi_1,\dots, \xi_e\}$ that contain no
subword from the set $\{\xi_e\xi_1,\dots,\xi_e\xi_r\}$ form a basis of
$\mathcal E$ over $k$.
   \item
Each bounded below (respectively, finitely presented) graded $\mathcal
E$-module has a resolution of length $2$ by bounded below (respectively,
finite) graded free modules.
   \end{enumerate}
    \end{theorem}

   \begin{proof}[Proof of Theorem {\rm \ref{thm:distinguished}}]
(1) and (2).  As $R$ is Koszul, the $k$-algebra $\mathcal E$ is  
generated
by $\mathcal E^1$, see \cite[(7.5)]{Se}, so there is a surjective
map $k\{\xi_1,\dots, \xi_e\}\to\mathcal E$ of graded $k$-algebras.
Sj\"odin, \cite[Thm.~4]{Sj1} shows that each $\phi_h$ is in its kernel,
so it induces a surjective homomorphism $\mathcal D\to \mathcal E$,  
where
we have set $\mathcal D=k\{\xi_1,\dots, \xi_e\}/(\phi_1,\dots,\phi_r)$.

Order the words in $\xi_1,\dots,\xi_e$ first by degree, and within each
degree lexicographically with respect to $\xi_e>\dots>\xi_1$.  This
order is admissible in the sense of \cite{Uf}.  By \cite[Thm.~7]{Uf},
a spanning set of $\mathcal D$ over $k$ is given by those words that
do not contain as a subword the leading monomial of any element from
the ideal $(\phi_1,\dots,\phi_r)$.  The leading term of $\phi_i$ is
$\xi_e\xi_i$ so $\mathcal D$ is spanned, \textit{a fortiori}, by the
words containing no subword $\xi_e\xi_i$ for $i=1,\dots,r$.  Call such
a word reduced. Let $w_n$ denote the number of reduced words of degree $n$; for $1\le
i\le e$, let $w_{n,i}$ denote the number of reduced words of degree $n$
ending in $\xi_i$. One then has
  \begin{alignat*}{4}
w_{n,i}&=w_{n-1} &&\text{for } &r+1&\le i\le e
  &\text{ and }n&\ge 1\,;
\\
  w_{n,i}& =w_{n-1}-w_{n-1,e}
   \quad &&\text{for } &1&\le i\le r &\text{ and }n&\ge 2\,.
  \end{alignat*}
For every $n\ge 2$ the definitions and the relations above yield
   \begin{align*}
w_n&=\sum_{i=1}^rw_{n,i}+\sum_{i=r+1}^ew_{n,i}\\
 &=r(w_{n-1}-w_{n-2})+(e-r)w_{n-1}\\
 &=ew_{n-1}-rw_{n-2}\,.
  \end{align*}
These equalities, along with $w_0=1$ and $w_1=e$, imply the second
equality below:
   \begin{align*}
\sum_{n=0}^\infty \rank_k(\mathcal E^n)t^n
  &=\frac 1{1-et+rt^2}\\
  &=\sum_{n=0}^\infty w_nt^n\\
  &\succcurlyeq\sum_{n=0}^\infty\rank_k(\mathcal D^n)t^n\\
  &\succcurlyeq\sum_{n=0}^\infty\rank_k(\mathcal E^n)s^n\,.
   \end{align*}
The first equality is \eqref{eq:koszul-field}, and the coefficient-wise
inequalities $\succcurlyeq$ of formal power series are evident.  Thus,
all the relations above are equalities.  It follows that the homomorphism
$\mathcal D\to \mathcal E$ is bijective and the reduced words form a basis.

(3)  The graded $k$-algebra $\agr R$ is Koszul, by the definition of the
Koszul property for $R$.  Thus, the $k$-algebra $\Ext_{\agr R}(k,k)$ is generated by
$\Ext_{\agr R}^1(k,k)$.  As $\fm^3=0$ holds, L\"ofwall
\cite[Thm.~2.3]{Lo} gives $\Ext_{\agr R}(k,k)\cong\mathcal E$ as graded
$k$-algebras.  Koszul duality yields $\Ext_{\mathcal E}(k,k)\cong\agr R$,
see e.g.\ \cite[4.1]{HI}, so $\Ext_{\mathcal E}^n(k,k)\cong\agr[n] R=0$
for $n>2$. This implies that each bounded below graded $\mathcal E$-module
has a resolution of length $2$ by bounded below graded free
modules; see \cite[\S8, Prop.~8, Cor.~5 and Cor.~2]{Bo}.

Theorem \ref{thm:factorization} yields a Golod homomorphism from a
hypersurface ring onto $R$.  By a result of Backelin and Roos, see
\cite[Thm.~4,~Cor.~3]{BR}, this implies that $\mathcal E$ is coherent.
It follows that when a graded $\mathcal E$-module $\mathcal M$ is finitely presented the kernel of the
map in any finite free presentation is a finite free $\mathcal E$-module.
  \end{proof}

Let $M$ be an $R$-module, $F$ a free resolution of $M$, and $G$ one of $k$.  Composition turns $\Hom_R(F,G)$ into a DG module over the DG algebra $\Hom_R(G,G)$. Thus, $\Ext_R(M,k)=\hh{\Hom_R(F,G)}$ is a graded left module over $\mathcal E=\hh{\Hom_R(G,G)}$ with $\Ext^0_R(M,k)=\Hom_R(M,k)$, and hence one has a natural homomorphism
\[
\mathcal E \otimes_k \Hom_R(M,k) \lra \Ext_R(M,k)
\]
of left $\mathcal E$-modules, which is an isomorphism when $\fm M=0$.

We index graded objects following custom and convenience:  For the $n$th component of a graded $k$-vector space $V$ we write either $V_n$ or $V^{-n}$.  Let $\shift^{i}V$ denote the graded vector space
with $(\shift^{i}V)_n=V_{n-i}$ for $n\in\BZ$; equivalently, $(\shift^{i}V)^n=V^{n+i}$. 

\begin{construction}
   \label{con:delta}
The $\agr R$-module structure on $\agr M$ defines a $k$-linear map
   \[
\nu\col \agr[1]R\otimes_k \agr[0]M\lra \agr[1]M
   \quad \text{with}\quad
\ov a\otimes \ov x\longmapsto \ov{ax}\,.
   \]
Set $-^*=\Hom_R(-,k)$.  As one has $\mathcal E^1=(\agr[1]R)^*$, one gets a $k$-linear map
   \[
(\agr[1]M)^*
  \xra{\nu^*\ }(\agr[1]R\otimes_k \agr[0]M)^*
  \cong (\agr[1]R)^*\otimes_k(\agr[0]M)^*=\mathcal E^1 \otimes_k(\agr 
[0]M)^*\,.
\]
As $\xi_1,\dots,\xi_e$ is a $k$-basis for $\mathcal E^1$, for each $ 
\psi\in
(\agr[1]M)^*$ there are uniquely defined elements $\psi_1,\dots,\psi_e$
in $(\agr[0]M)^*$ such that $\nu^*(\psi)= \sum_{i=1}^e \xi_i \otimes
\psi_i$. Evidently, the map
   \[
\delta \col \mathcal E\otimes_k \shift^{-1}(\agr[1]M)^*\lra \mathcal
E\otimes_k (\agr[0]M)^* \quad\text{with}\quad \delta (\xi \otimes
\psi)= \sum_{i=1}^e \xi\xi_i \otimes  \psi_i
   \]
is an $\mathcal E$-linear homomorphism of degree zero.
     \end{construction}

\begin{theorem}
  \label{thm:distinguished-modules}
Let $(R,\fm,k)$ be a local ring such that $\fm$ has a Conca generator,
let $M$ be a finite $R$-module with $\fm^2M=0$, and set 
\[
\mathcal F=\Ker (\delta)\,.
\]

The graded $\mathcal E$-module $\mathcal F$ then is finite free, it
satisfies
   \begin{align}
   \label{eq:bound}
\mathcal F^n&=0
  \quad\text{for}\quad
 n\le1\,;
   \\
 \min\{j\ge 0\mid (k\otimes_{\mathcal E}\mathcal F)^{\ges j+2}=0\}
   &=\min\{n\ge0\mid\Omega^R_n(M)\text{ is Koszul\,}\}\,,
  \end{align}
and the following sequence is a minimal free resolution of
$\Ext_R(M,k)$ over $\mathcal E$:
   \begin{equation}
   \label{eq:resolution}
\xymatrixcolsep{1pc} \xymatrix{
     0\ar@{->}[r]
&\mathcal F\ar@{->}[r]
&\mathcal E\otimes_k \shift^{-1}(\agr[1]M)^*
\ar@{->}[rr]^-{\left[\begin{smallmatrix}\delta \\ 0\end{smallmatrix} 
\right]}
&&\mathcal E\otimes_k (\agr[0]M)^*\oplus\shift^{1} \mathcal F\ar@{->}[r]
&0
   }
   \end{equation}
      \end{theorem}

   \begin{remark}
  \label{rmk:remove-hyp}
The proof shows that the conclusions above, except for the
finiteness of $\mathcal F$ over $\mathcal E$, hold when the hypothesis
that $R$ has a Conca generator is weakened to assuming that $R$ is  
Koszul
and $\fm^3=0$; see also Roos \cite[(3.1)]{Ro} for \eqref{eq:resolution}.
   \end{remark}

\begin{proof}[Proof of Theorem~\emph{\ref{thm:distinguished-modules}}]
Since ${\mathcal E}^{n}=0$ for $n<0$, it follows from the definition of $\delta$ that ${\mathcal F}^{n}=0$ for $n<0$; see \ref{con:delta}. Moreover, $\delta^{0}$ is the map $\nu^{*}\col {\agr[1]M}^{*}\to {\agr[1]R}^{*}\otimes_{k}{\agr[0]M}^{*}$, which is injective, since $\nu$ is surjective. Thus ${\mathcal F}^{0}=0$ as well.
Thus \eqref{eq:bound} is proved.

Set $\mathcal M=\Ext_R(M,k)$. The exact sequence of $R$-modules
   \[
0\lra \agr[1]M\lra M\lra \agr[0]M\lra0
   \]
induces an exact sequence of graded left $\mathcal E$-modules
   \[
   \xymatrixrowsep{.7pc} 
   \xymatrixcolsep{1.5pc} \xymatrix{
\shift^{-1} \Ext_R(\agr[1]M,k)\ar@{->}[r]^-{\eth}
   &\Ext_R(\agr[0]M,k)\ar@{->}[r]
   &\mathcal M\ar@{->}[r]
   &
   \\
{\hphantom{xxx}}\Ext_R(\agr[1]M,k)\ar@{->}[r]^-{\shift\eth}
   &\shift^{1}\Ext_R(\agr[0]M,k)
   }
   \]
For $j=0,1$ one has $\Ext_R(\agr[j]M,k)\cong {\mathcal
E}\otimes_k(\agr[j]M)^*$ as graded left $\mathcal E$-modules.

Using Proposition \ref{diagram} below, it is not hard to check
that one can replace the map $\eth$ with $\delta$.  Theorem
\ref{thm:distinguished}(3) then implies that $\mathcal F$ is finite
projective over $\mathcal E$.  Bounded below projective graded $\mathcal
E$-modules are free, see \cite[\S8, Prop.~8, Cor.~1]{Bo}, so the exact
sequence above yields the free resolution \eqref{eq:resolution}.  From this
resolution one gets
   \[ \max\{j\ge0\mid\Tor^{\mathcal E}_n(k,\mathcal M)_{n+j}\ne0\text{
   for some }n\} =\min\{j\ge 0\mid (k\otimes_{\mathcal E}\mathcal F)^ 
{\ges
   j+2}=0\}\,.
  \]
By \cite[(5.4)]{HI}, the number on the left-hand side of the equality above
is equal to the least integer $n\ge0$, for which the $R$-module
$\Omega^R_n(M)$ is Koszul, so we are done.
      \end{proof}

\begin{corollary}
  \label{cor:distinguished-modules}
For every finite $R$-module $L$ the graded $\mathcal E$-module
$\Ext_R(L,k)$ has a finite free resolution of length at most $2$.
      \end{corollary}

      \begin{proof}
Let $R^m\to L$ be a free cover.  The exact sequence of $R$-modules
   \[
0\lra K\lra R^m\lra L\lra0
   \]
induces an exact sequence of graded $\mathcal E$-modules
     \[
  0\lra \shift^{-1}\Ext_R(K,k)\lra\Ext_R(L,k)\lra k^m\lra0
       \]
One has $\fm^2K\subseteq\fm^2(\fm R^m)= \fm^3R^m=0$, so the graded
$\mathcal E$-module $\Ext_R(K,k)$ is finitely presented by the theorem.
The ideal $\mathcal E^{\ges1}$ of $\mathcal E$ is finitely generated,
by Theorem \ref{thm:distinguished}(3), so the $\mathcal  E$-module $k 
$ is
finitely presented, and hence so is  $k^m$.  It follows that $\Ext_R 
(L,k)$
is finitely presented.  Now refer to Theorem \ref{thm:distinguished}(3).
      \end{proof}

We have deferred to the end of the section the statement and proof of
a simple general homological fact, for which we could not find an
adequate reference.

\begin{proposition}
\label{diagram}
Let $R$ be a commutative ring, $M$ an $R$-module, $I$ an ideal,
and $\nu\col R_1\otimes_{R_0}M_0\to M_1$ the natural map, where
$R_j=I^j/I^{j+1}$ and $M_j=I^jM/I^{j+1}M$.

If $M_0$ is free over $R_0$, then there is a commutative diagram of
$R_0$-linear maps
\[
\xymatrixrowsep{2pc}\xymatrixcolsep{5pc} \xymatrix{
\Hom_R(M_1,{R_0})
   \ar@{->}[r]^-{\eth}
    \ar@{->}[d]_{\cong}
&\Ext^1_R(M_0,{R_0})
    \ar@{->}[d]^{\cong}
\\
\Hom_{R_0}(M_1,{R_0})
    \ar@{->}[r]^-{\Hom_{R_0}({\nu},{R_0})}
&\Hom_{R_0}(R_1\otimes_{R_0} M_0,{R_0})
}
   \]
where $\eth$ is the connecting homomorphism defined by the exact
sequence
   \[
0\lra M_1\lra M/I^2M\lra M_0\lra 0
   \]
     \end{proposition}

   \begin{proof}
The hypothesis allows us to choose a surjective homomorphism of
$R$-modules $F\to M_0$ with $F$ free, so that the induced map $F/IF\to
M_0$ is an isomorphism.

The exact sequence above appears in a commutative diagram with exact  
rows
   \[
\xymatrixrowsep{1.1pc} \xymatrixcolsep{2.8pc} \xymatrix{
0
   \ar@{->}[r]
&M_1
   \ar@{->}[r]
& M/I^2M
   \ar@{->}[r]
&M_0
   \ar@{->}[r]
&0
   \\
&R_1\otimes_{R_0}M_0
   \ar@{->}[u]^-{\nu}
   \\
0
   \ar@{->}[r]
&I\otimes_RF
   \ar@{->}[r]
   \ar@{->}[u]
& F
   \ar@{->}[r]
   \ar@{->}[uu]
&  M_0
   \ar@{->}[r]
   \ar@{=}[uu]
&0
}
\]
where the bottom row is the result of tensoring with $F$ the
exact sequence of $R$-modules $0\to I\to R\to R_0\to0$.  The isomorphisms in the  
induced
commutative diagram
    \[ \xymatrixrowsep{1.1pc}\xymatrixcolsep{.76pc}
\xymatrix{ &&\Hom_R(M_1,{R_0})
   \ar@{->}[r]^-{\eth} \ar@{->}[d]_-{\Hom_{R}(\nu,{R_0})}
&\Ext^1_R(M_0,{R_0})
    \ar@{=}[dd]
\\ &&\Hom_{R}(R_1\otimes_{R_0} M_0,{R_0})
   \ar@{->}[d]^-{\cong}
\\ \Hom_R(M_0,{R_0})
    \ar@{->}[r]^-{\cong}
&\Hom_R(F,{R_0})
   \ar@{->}[r]
&\Hom_R(I \otimes_R F,{R_0})
   \ar@{->}[r]^-{\eth'}
&\Ext^1_R(M_0,{R_0})
   \ar@{->}[r]
&0
  }
   \]
are due to the isomorphism $F/IF\xra{\cong} M_0$.  The
exactness of the bottom row shows that the connecting map
$\eth'$ is bijective.  Furthermore, the following diagram
\[ \xymatrixrowsep{2pc}\xymatrixcolsep{2.5pc} \xymatrix{
\Hom_{R_0}(M_1,{R_0})
   \ar@{->}[d]_-{\Hom_{R_0}(\nu,{R_0})} \ar@{->}[r]^-{\cong}
&\Hom_{R}(M_1,{R_0})
   \ar@{->}[d]^-{\Hom_{R}(\nu,{R_0})}
\\ \Hom_{R_0}(R_1\otimes_{R_0} M_0,{R_0})
   \ar@{->}[r]^-{\cong}
&\Hom_{R}(R_1\otimes_{R_0} M_0,{R_0})
  }
   \]
commutes by functoriality. The horizontal arrows are isomorphisms
because $R\to R_0$ is surjective.  
One gets the desired result by combining the last two diagrams 
   \end{proof}

\section{Koszul modules}
\label{sec:Koszul modules}

Theorem \ref{thm:koszuls} shows that when $\fm$ has a Conca generator
the asymptotic properties of arbitrary free resolutions are determined by
those of resolutions of Koszul modules.  In this section we turn to the
problem of identifying and exhibiting such modules.  The next result
follows easily from work in the preceding section.

\begin{proposition}
  \label{prop:koszul}
Let $(R,\fm,k)$ be a Koszul local ring with $\fm^3=0$. Set
\[
e=\rank_k(\fm/\fm^2)\,,\quad r=\rank_k(\fm^2)\,,\quad\text{and}\quad
\mathcal E=\Ext_R(k,k)\,.
\]
For each finite $R$-module $M$ with $\fm^2M=0$ the following are
equivalent.
   \begin{enumerate}[\quad\rm(i)]
  \item
The $R$-module $M$ is Koszul.
  \item
The Poincar\'e series of $M$ over $R$ is given by the formula
   \[
\displaystyle\po{M}t=\frac{\hilb M{-t}}{1-et+rt^2}\,.
  \]
  \item
The $\mathcal E$-module $\Ext_R(M,k)$ has projective dimension at
most $1$.
   \item
The map $\delta$ from Construction \emph{\ref{con:delta}} is injective.
   \end{enumerate}
      \end{proposition}

  \begin{proof}
By Remark \ref{rmk:remove-hyp}, the resolution \eqref{eq:resolution}
yields (iv) $\iff$ (iii), and an equality
  \[
\po{M}t=\frac{\hilb M{-t}}{1-et+rt^2}
+\left(1+\frac1t\right)\cdot \sum_{n=0}^\infty\rank_k(\mathcal F^n)\,t^n 
\,,
   \]
from which (iv) $\iff$ (ii) follows.  The equality \eqref{eq:bound}
establishes (iv) $\iff$ (i).
   \end{proof}

Next we exhibit a substantial family of Koszul modules.

\begin{theorem}
   \label{thm:annihilated}
If $(R,\fm,k)$ is a local ring and $x$ is a Conca generator of $\fm$,
then every finite $R$-module $M$ annihilated by $x$ is Koszul.
   \end{theorem}

For the proof we need a general change-of-rings result for the Koszul
property.

\begin{lemma}
  \label{lem:generic-lifting}
  Let $(R,\fm,k)$ and $(R',\fm',k)$ be Koszul local rings, $\rho\col  
R'\to R$ a Golod
  homomorphism, and $M$ a finite $R$-module.

If $M$ is Koszul over $R'$, then it is Koszul over $R$ as well.
   \end{lemma}

   \begin{proof}
In view of \cite[(6.1)]{HI}, it suffices to show that the module $M$
is $\rho$-Golod in the sense of Levin \cite[(1.1)]{Lv}; that is, the
map $\tor[n]{\rho}Mk$ is injective for each $n\ge0$.

Fix a non-negative integer $n$.  The exact sequence of $R$-modules
   \[
     0\lra \fm M\xrightarrow{\,\iota\,} M\xrightarrow{\,\pi\,} M/\fm M 
\lra 0
   \]
induces an exact sequence of $k$-vector spaces
   \[
     \Tor^{R'}_n(\fm M,k)\xrightarrow{\Tor^{R'}_n(\iota,k)}
     \Tor^{R'}_n(M,k)\xrightarrow{\Tor^{R'}_n(\pi,k)}\Tor^{R'}_n(M/ \fm M,k)\,.
   \]
When $M$ is Koszul over $R'$ one has $\Tor^{R'}_n(\iota,k)=0$ by
\cite[(3.2)]{Se}, and hence the map $\Tor^{R'}_n(\pi,k)$  is injective.
It appears in a commutative diagram
   \[
  \xymatrixrowsep{2.3pc} \xymatrixcolsep{1.5pc}
   \xymatrix{
\tor[n]{R'}{M}k
\ar@{->}[rr]^-{\tor[n]{R'}{\pi}k}
\ar@{->}[d]_-{\tor[n]{\rho}{M}k}
  &&
\tor[n]{R'}{M/\fm M}k
\ar@{->}[r]^-{\cong}
\ar@{->}[d]^-{\tor[n]{\rho}{M/\fm M}k}
   &
(M/\fm M)\otimes_k\tor[n]{R'}kk
\ar@{->}[d]^-{(M/\fm M)\otimes_k\tor[n]{\rho}kk}
    \\
\tor[n]{R}{M}k
\ar@{->}[rr]^-{\tor[n]{R}{\pi}k}
&&
\tor[n]{R}{M/\fm M}k
\ar@{->}[r]^-{\cong}
&
(M/\fm M)\otimes_k\tor[n]{R}kk
}
   \]
Since $\rho$ is Golod, $\tor[n]{\rho}kk$ is injective by
\cite[(3.5)]{A0}, hence so is $\tor[n]{\rho}{M}k$.
   \end{proof}

The following explicit construction is also used in the  proof of the
theorem.

\begin{remark}
   \label{rem:pol-presentation}
Choose a minimal generating set $\{x_{1},\dots,x_{e}\}$ of $\fm$
as in \ref{rem:quadrics}, and a map $\pi\col(P,\fp,k)\to R$ as in
Theorem \ref{thm:factorization}.  Pick elements $u_{1},\dots,u_{e}$
in $\fp$ so that $\pi(u_i)=x_i$, and let $X_i$ denote the leading
form of $u_i$ in the associated graded ring $\agr P$.  Thus, one has
$\agr P=k[X_1,\dots,X_{e}]$, and the elements $X_1,\dots,X_{e}$ are
algebraically independent over $k$.  Let $I$ be the ideal of $\agr P$
generated by
   \begin{align}
   \label{eq:quadrics}
X_iX_j-\sum _{h=1}^r a^{ij}_hX_hX_e
   &\quad\text{for}\quad 1\le i\le j\le e-1\,,
\quad\text{and}
   \\
   \label{eq:monomials}
X_{l}X_e
   &\quad\text{for}\quad r+1\le l\le e
   \end{align}
where $r=\rank_k\fm^2$ and the elements $a^{ij}_h\in k$ are defined by
formula \eqref{eq:quadrics-local}.

Let $\ov u_i$ denote the image of $u_i$ in the ring $\agr P/I$.  As 
$\ov u_e$ is a Conca generator of $(\ov u_1,\dots,\ov u_{e})$, one 
has $(\ov u_1,\dots,\ov u_{e})^3=0$, hence $\rank_k(\agr P/I)=1+e+r$.
Since $I$ is contained in the kernel of the map $\agr\pi\col\agr P
\to\agr R$ of graded $k$-algebras, it induces a surjection $\agr P/I\to 
\agr R$, which is bijective because $\rank_k\agr R=1+e+r$ holds.
  \end{remark}

\begin{proof}[Proof of Theorem {\rm\ref{thm:annihilated}}]
Let $M$ be a finite $R$-module with $xM=0$.

To prove that $M$ is Koszul it suffices to show that the graded $\agr
R$-module $\agr M$ has a linear free resolution; see \cite[(2.3)]{Se} or
\cite[(1.5)]{HI}.  By Remark \ref{rem:pol-presentation}, the ring
$\agr R$ is local and the initial form $\ov x$ of $x$ is
a Conca generator of its maximal ideal.  As $\ov x\agr M=0$ evidently
holds, after changing notation we may assume that $R$ is graded, $\fm$
has a Conca generator $x\in R_1$, $M_j=0$ for $j\ne0,1$, and  
$M_1=R_1M_0$.

Remark \ref{rem:pol-presentation} yields an isomorphism $R\cong
k[X_1\dots,X_e]/I$, where $I$ is generated by the quadratic forms
in \eqref{eq:quadrics} and \eqref{eq:monomials}.  Thus, there is a
surjective homomorphism $\rho\col R'\to R$ of graded $k$-algebras,
where $R'=k[X_1\dots,X_e]/I'$ and $I'$ is the ideal generated by $X_e^2$
and the forms in \eqref{eq:quadrics}.  The image of $X_e$ is a
Conca generator of the ideal $\fm'=(X_1\dots,X_e)R'$; in particular,
$R'$ is local with maximal ideal $\fm'$.

Theorem \ref{thm:koszuls} shows that both rings $R$ and $R'$ are Koszul.
Furthermore, one has an isomorphism $\Ker(\rho)\cong k(-2)^{e-r-1}$ of
graded $R'$-modules, so $\Ker(\rho)$ has a $2$-linear free resolution
over $R'$.  This implies that the homomorphism $\rho$ is Golod,
see \cite[(5.8)]{HI}.  Referring to Lemma \ref{lem:generic-lifting}
one sees that to finish the proof of the theorem it suffices to show
that $M$ has a linear free resolution over $R'$.  Replacing $R$ with $R'$
and changing notation once more, we may also assume $\rank_k(R_2)=e-1$.

Next we prove $xR=(0:x)_R$.  The condition $x^2=0$ implies
$xR\subseteq(0:x)_R$.  Equality holds because both ideals have the same
rank: The exact sequences
   \begin{gather*}
0\lra \fm^2\lra xR\lra xR/\fm^2\lra 0
   \\
0\lra (0:x)_R\lra R\lra xR\lra 0
   \end{gather*}
yield $\rank_k xR=(e-1)+1=e$ and $\rank_k(0:x)_R=2e-e=e$.  The  
complex
   \[
\cdots\lra
R(-2)\xrightarrow{\,x\,}R(-1)\xrightarrow{\,x\,}R\lra0\lra\cdots
   \]
is thus a minimal free resolution of the graded $R$-module $R/xR$, and
hence $R/xR$ is Koszul.  As one has $xM=0$, there is an exact sequence
of graded $R$-modules
   \[
0\lra k^q(-1)\lra (R/xR)^m\lra M\lra0
   \]
with $m=\rank_k (M/\fm M)$.  It induces for each pair $(n,j)$
an exact sequence
   \[
\gtor n{R}{(R/xR)^m}kj\lra\gtor n{R}{ M}kj\lra\gtor{n-1}{R}{k^q}k{j-1}\,.
   \]
The vector spaces at both ends are zero when $n\ne j$, because both $k$
and $R/xR$ are Koszul over $R$.  As a consequence, we get $\gtor n{R}{
M}kj=0$ for $n\ne j$, as desired.
   \end{proof}

\section{Gorenstein rings and related rings}
\label{sec:Gorenstein rings and related rings}

   \setcounter{equation}{0}

In this section we study modules over local rings $(R,\fm,k)$ with
$\fm^3=0$ and $\rank_k(\fm^2)\le1$, focusing on the Koszul property.
We prove that, outside of an easily understood special case, every
module has a Koszul syzygy.   We completely describe the non-Koszul
indecomposable modules when $R$ is Gorenstein.  An important finding is
that for $e=\rank_k(\fm/\fm^2)$ the sequence $(b_n)_{n\ges0}$ defined by  
   \begin{equation}
      \label{eq:recursion}
b_n=
   \begin{cases}
n+1
  &\text{when }e=2\,;\\
\displaystyle\frac1{2^{n+1}}\sum_{j=0}^{\lfloor{(n-1)}/2\rfloor} 
(e^2-4)^{j}e^{n-2j}
  &\text{when }e\ge3\,,
   \end{cases}
     \end{equation}
provides numerical invariants for checking the  Koszul  property of
$R$-modules.

\begin{theorem}
   \label{thm:r1}
Let $(R,\fm,k)$ be a local ring with $\fm^3=0$ and $\rank_k\fm^2=1$.

For $e=\rank_k(\fm/\fm^2)$ and $s=\rank_k(0:\fm)$ the following
are equivalent.
   \begin{enumerate}[\quad\rm(i)]
  \item
There is an inflation $R\to(R',\fm',k')$ such that $\fm'$ has a Conca
generator.
  \item
$R$ is Koszul.
  \item
$\po kt\cdot(1-et+t^2)=1$.
   \item
$\beta^R_n(k)=b_n$ for every integer $n\ge0$.
  \item
$s\le e-1$.
   \end{enumerate}
    \end{theorem}

The most useful property above is (i), as it has consequences for
free resolutions of all $R$-modules; see Theorem \ref{thm:koszuls} and
Corollary \ref{cor:sjodin}.  Thus, the main thrust of the theorem
is that (i) follows from the easily verifiable condition (v).

The equivalence of conditions (ii), (iii), and a modified form of (v)
is established by  Fitzgerald in \cite[(4.1)]{Fi}.  From the proof of
that result we abstract the following statement, which it is not hard
to verify directly.

   \begin{remark}
    \label{rem:fitz}
A local ring $(R,\fm,k)$ with $\fm^3=0$ and $\rank_k\fm^2=1$ is isomorphic to a fiber
product $S\times_kT$, where $(S,\fs,k)$ is a local ring with $\fs^2=0$
and $(T,\ft,k)$ is a Gorenstein local ring with $\ft^2\ne0$; see the
proof that (1) implies (3) in \cite[(4.1)]{Fi}.

Set $p=\rank_k(\fs/\fs^2)$ and $q=\rank_k(\ft/\ft^2)$ and note the  
evident
relations
   \[
\rank_k(\fm/\fm^2)=p+q
   \quad\text{and}\quad
\rank_k(0:\fm)=p+1\,.
   \]
   \end{remark}

We also note an elementary observation, to be used more than once.

   \begin{remark}
    \label{rem:socle}
A Gorenstein local ring $(R,\fm,k)$ with $\fm^3=0\ne\fm^2$ satisfies
   \[
k\cong (0:\fm)=\fm^2 = x\fm
   \quad\text{for every}\quad
x\in\fm\smallsetminus\fm^2\,.
  \]

Indeed, the Gorenstein property of $R$ implies $k\cong(0:\fm)$.  As the
hypotheses on $\fm$ mean $(0:\fm)\supseteq\fm^2\ne0$, one concludes
$(0:\fm)=\fm^2$.  Thus, one has either $x\fm =\fm^2$ or $x\fm =0$,
but the second option entails $x\in(0:\fm)=\fm^2$, a contradiction.
   \end{remark}

\begin{proof}[Proof of Theorem \emph{\ref{thm:r1}}]
(i) $\implies$ (ii).  This holds by Theorem \ref{thm:koszuls}.

(ii) $\implies$ (iii). This follows from Remark \ref{rem:syzygy}.

(iii) $\iff$ (iv).  This is seen by decomposition into prime fractions.

In the rest of the proof we use the notation of Remark \ref{rem:fitz}.

(iii) $\implies$ (v).  The hypothesis and Remark \ref{rem:dress} yields
equalities
   \[
1-et+t^2=\frac1{\po kt}= \frac1{\po[T]kt} + \frac1{\po[S]kt}-1 = \frac1{\po[T]kt}-pt\,.
   \]
They imply $\po[T]kt^{-1}=1-qt+t^2$. This rules out the case $q=1$, because 
the local ring $T$ then has embedding dimension $1$ so one has $\po[T]kt^{-1}=1-t$. 
Thus $s=p+1= e-q+1\le e-1$.

(v) $\implies$ (i).
There exists a local ring $(R',\fm',k')$ with $k'$ algebraically closed
and an inflation $R\to R'$, see \cite[Ch.~0, 10.3.1]{Gr}.  One then has
$\fm'^3=0$ and $\rank_{k'}\fm'^2=1$; also, $\rank_{k'}(\fm'/\fm'^2)=e$
and $\rank_{k'}(0:\fm')=s$ hold.  Thus, it suffices to prove that
$\fm$ has a Conca generator when $k$ is algebraically closed.

The Gorenstein ring $T$ in the decomposition $R=S\times_kT$ has
$q=e-p=e-(s-1)\ge2$.   A Conca generator of $\ft$ clearly also is a  
Conca
generator of $\fm$.  Thus, we may further assume that $R$ is Gorenstein
with $e\ge2$; this implies $\fm^2\ne0$.

Remark \ref{rem:socle} shows that for each $x\in\fm\smallsetminus\fm^2$
one has $x\fm=\fm^2$.  On the other hand, by \cite[Lem.~3]{Co1} one can
choose an element $x$ as above, so that $\ov x\in\fm/\fm^2=\agr[1]{R}$
satisfies $\ov x^2=0$.  This yields $x^2\in\fm^3=0$, so $x$ is  a Conca
generator of $\fm$.
    \end{proof}

Next we recover and extend Sj\"odin's description \cite{Sj2} of  
Poincar\'e series of modules over Gorenstein local rings with $\fm^3=0$.  Note 
that the second case below differs from the other two, as all series
have a common denominator \emph{different from} $\hilb Rt$.

   \begin{corollary}
  \label{cor:sjodin}
Let $(R,\fm,k)$ be a local ring with $\fm^3=0$ and $\rank_k\fm^2\le1$.

The numbers $e=\rank_k(\fm/\fm^2)$ and $s=\rank_k(0:\fm)$ then satisfy
$s\le e$, and for every finite $R$-module $M$ one has
     \begin{alignat*}{2}
\deg\big(\po Mt\cdot(1-et)\big)&\le1\ \  &&\text{ when $\fm^2=0$.}
   \\
\deg\big(\po Mt\cdot(1-et)\big)&\le2 &&\text{ when $\rank_k\fm^2=1$
and $s=e$.}
   \\
\deg\big(\po Mt\cdot(1-et+t^2)\big)&<\infty &&\text{ when $\rank_k 
\fm^2=1$
and $s<e$.}
     \end{alignat*}
   \end{corollary}

   \begin{proof}
When $\fm^2=0$ one has $\Omega^R_1(M)\cong k^a$ for some $a\ge0$, and
hence
   \[
\po Mt-\beta^R_0(M)=\po{k^a}t\cdot t=\frac{a}{1-et}\cdot t
   \]

When $\rank_k\fm^2=1$ and $s=e$ Remark \ref{rem:fitz}
yields $R=S\times_kT$, where the local ring $(T,\ft,k)$ has
$\rank_k(\ft/\ft^2)=1$.  Thus, the maximal ideal $\ft$ is principal,
so every indecomposable $T$-module is isomorphic to $k$, $\ft$, or $T$.

By \cite[Rem.~3]{DK} one has $\Omega^R_2(M)=K\oplus L$, where $K$ is
an $S$-module and $L$ is a $T$-module.  The inclusion $\Omega^R_2(M)
\subseteq\fm F_1$, where $F_1$ is a free $R$-module, gives $\fs
K\subseteq \fs^2 F_1=0$ and $\ft^2 L\subseteq \ft^3 F_1=0$.  Thus,
there exist integers $a,b,c\ge0$ and isomorphisms $K\cong k^a$ and
$L\cong k^b\oplus\ft^c$ of $S$-modules and $T$-modules, respectively.

For $d=e-1$ the discussion above leads to the first and third equalities
below:
   \begin{align*}
\po Mt-\beta^R_0(M)-\beta^R_1(M)\cdot t
&=\po Kt\cdot t^2+\po Lt\cdot t^2
   \\
&=\frac{\po[S]{K}t\cdot\po[T]kt+\po[S]kt\cdot\po[T]Lt}
{\po[S]kt+\po[T]kt-\po[S]kt\cdot\po[T]kt}\cdot t^2
   \\
&=\frac{\displaystyle\frac
a{1-dt}\cdot\frac 1{1-t}+\frac 1{1-dt}\cdot\frac{b+c}{1-t}}
{\displaystyle\frac 1{1-dt}+\frac 1{1-t}-\frac 1{1-dt}\cdot\frac
1{1-t}}\cdot t^2
   \\
&=\frac{a+b+c}{1-et}\cdot t^2
   \end{align*}
The second equality comes from the change-of-rings formula in Remark
\ref{rem:dress}.

When $\rank_k\fm^2=1$ and $s<e$ the result follows from Theorems
\ref{thm:r1} and \ref{thm:koszuls}.
    \end{proof}

It follows from work of Conca that, under additional hypotheses, certain
conditions of Theorem \ref{thm:r1} remain valid for larger values of
$\rank_k\fm^2$.

    \begin{remark}
Let $(R,\fm,k)$ be a local ring with $\fm^3=0$, field $k$ of  
characteristic
different from $2$, and $k$-algebra $\agr R$ defined by quadratic  
relations.

When $\rank_k\fm^2=2$ holds there is an inflation $R\to(R',\fm',k')$,
such that $\fm'$ has a Conca generator; when $\rank_k\fm^2=3$ holds the
ring $R$ is Koszul.

Indeed, as in the proof of (v) $\implies$ (i) in Theorem \ref{thm:r1}  
one
may reduce to the case when $R=\agr R$ and $k$ is algebraically closed.
If $\rank_kR_2=2$, then the proof of \cite[Prop.~6]{Co1} shows that the
ideal $(R_1)$ has a Conca generator.  If $\rank_kR_2=3$, then the ring
$R$ is G-quadratic by \cite[Thm. 1.1(2)]{Co2}; in particular, $R$ is a Koszul ring.
    \end{remark}

We turn to Koszul modules over a Gorenstein local ring $(R,\fm,k)$
with $\fm^3=0$.

Recall that when $\rank_k(\fm/\fm^2)=1$ the ring $R$ has finite representation
type:  When $\fm^2=0$ the indecomposable $R$-modules are $k$ and $R$, and
both are Koszul.  When $\fm^2\ne0$ the indecomposable modules are $k$, $\fm$,
and $R$; only $R$ is Koszul.

When $\rank_k(\fm/\fm^2)=2$, the field $k$ is algebraically closed, and
$\operatorname{char}(k)R=0$, one has $R\cong k[X_1,X_2]/(X_1^2,X_2^2)$.
The ring $R$ has tame representation type, and its indecomposable
modules are described by Kronecker's classification of pairs of  
commuting matrices.  From this description one can deduce that the 
negative  syzygies of $k$ are the only non-Koszul indecomposable 
$R$-modules.

We prove that the latter property holds for all Gorenstein local rings
with $\fm^3=0$ and $e\ge2$; this may surprise, as their representation
theory is wild when $e\ge3$.

\begin{theorem}
   \label{thm:gorenstein}
Let $(R,\fm,k)$ be a Gorenstein local ring with $\fm^3=0$, and set
   \[
\Omega^R_{-n}(k)=\Hom_R(\Omega_n^R(k),R)
   \quad\text{for each}\quad
n\ge1\,.
   \]
Set $e=\rank_k(\fm/\fm^2)$, assume $e\ge2$ holds, and for $n\ge0$
define $b_n$ by \eqref{eq:recursion}.

A finite $R$-module $M$ is Koszul if and only if it has no direct  
summand
isomorphic to $\Omega^R_{-i}(k)$ with $i\ge1$; when $M$ is  
indecomposable
the following are equivalent.
   \begin{enumerate}[\quad\rm(i)]
  \item
$M$ is not Koszul.
   \item
$M\cong\Omega_{-i}^R(k)$ for some $i\ge1$.
   \item
$\hilb Mt=b_{i-1}+b_{i}\,t$.
   \end{enumerate}
\end{theorem}

\begin{proof}
It suffices to verify the equivalence of conditions (i) through (iii)  when
$M$ is an indecomposable, non-free, finite $R$-module.  These 
properties imply $\fm^2M=0$.

Indeed, by Remark \ref{rem:socle} one has $\fm^2=(0:\fm)\cong k$.  
Let $s$ be a generator of $\fm^2$ and $y\in M$ an element
with $sy\ne0$.  The $R$-linear map $R\to M$ sending $1$ to $y$ is
injective since it is injective on $(0:\fm)$. It is then split, for $R$
is self-injective, and so $M$ has a direct summand isomorphic to $R$;
a contradiction.

(i) $\implies$ (ii).
Proposition \ref{prop:koszul} yields $\po Mt\ne\hilb M{-t}/\hilb R{-t}$.  
By a result of  Lescot, see \cite[3.4(1)]{Le}, then $k$ is isomorphic to a direct 
summand of $\Omega_i^R(M)$ for some $i\ge1$.  As $\Omega_i^R(M)$ is 
indecomposable along with $M$, see \cite[1.3]{He}, one gets 
$\Omega_i^R(M)\cong k$, hence $M\cong\Omega^R_{-i}(k)$ by Matlis duality.

(ii) $\implies$ (iii).
A minimal free resolution $F$ of $M$ yields an exact sequence
     \begin{equation}
       \label{eq:syzygy-k}
0\lra k\lra F_{i-1}\lra F_{i-2}\lra\cdots\lra F_0\lra M\lra 0\,.
   \end{equation}
Set $G_n=\Hom_R(F_{i-1-n},R)$.  Applying $\Hom_R(-,R)$ one gets 
an exact sequence
     \begin{equation*}
0\lra \Hom_R(M,R)\lra G_{i-1}\lra G_{i-2}\lra\cdots\lra G_0\lra k\lra 0
   \end{equation*}
with $\dd(G_n)\subseteq\fm G_{n-1}$ for $n=0,\dots,i-1$.  From it
and Theorem \ref{thm:r1} one obtains
   \[
\rank_k(M/\fm M)=\rank_R(F_0)=\rank_R(G_{i-1})=b_{i-1}\,.
   \]
Note that $\fm M\ne 0$; else $b_i=1$ for some $i\ge 1$, which cannot
be the case. 

Since $\fm^2M=0$, one has that $\fm M\subseteq (0:\fm)_{M}$. The reverse inclusion holds because
the composed map $(0:\fm)_{M}\to M\to M/\fm M$ has to be zero, as $M$ is indecomposable.
Thus, $\fm M=(0:\fm)_M$, which gives the first equality below; Matlis duality
gives the second one, and \eqref{eq:syzygy-k} the third:
   \[
\rank_k(\fm M)=\rank_k((0:\fm)_M)=\rank_k(k\otimes_R\Omega^R_{i}(k))
=\rank_R(G_{i})=b_{i}\,.
   \]

(iii) $\implies$ (i).
Assuming that $M$ is Koszul, from Proposition \ref{prop:koszul} one gets
   \[
   \po Mt=(b_{i-1}-b_{i}t)\cdot\po kt=
   (b_{i-1}-b_{i}t)\cdot\sum_{n=0}^\infty b_nt^n\,,
   \]
hence $\beta_i^R(M)=b_{i-1}b_i-b_ib_{i-1}=0$.  Thus, $M$ has finite
projective dimension.  It is free because $R$ is artinian, contradicting
the hypotheses $\fm^2M=0\ne M$.
   \end{proof}

  \begin{corollary}
Let $(R,\fm,k)$ be a Gorenstein local ring with $\fm^3=0$. 

Set $e=\rank_{k}(\fm/\fm^{2})$. If $M$ is an indecomposable finite $R$-module such that
an inequality $\rank_k(\fm M)\le (e-1)\rank_k(M/\fm M)$ holds, then $M$ is Koszul.
  \end{corollary}

   \begin{proof}
The equivalence of (iii) and (iv) in Theorem \ref{thm:r1}
means that $(b_n)_{n\ge0}$ satisfies the recurrence relation
$b_{n+1}=eb_n-b_{n-1}$ for $n\ge2$, with $b_{0}=1$ and $b_{1}=e$.
It implies an inequality $b_{n}>(e-1)b_{n-1}$ for each $n\ge0$, so $M$ fails test
(iii) of Theorem \ref{thm:gorenstein}.
   \end{proof}

\section*{Acknowledgments}
   We thank Aldo Conca for useful remarks regarding this work, and the referee for a careful reading of the manuscript.

  \end{document}